\documentclass[12pt]{article}
\usepackage{amssymb}
\usepackage{amscd}
\usepackage{amsmath}
\usepackage{amstext}
\usepackage[all]{xy}

\title{\Large \bf  Balance with Unbounded Complexes
\thanks{2010 {\it Mathematics Subject Classification}.55U15,16E05,16E30,
18G15.}
\thanks {{\it Keywords}. Complexes, balance, Tate (co)homology.}}

\author{Edgar E. Enochs$^{1}$, Sergio Estrada$^{2}$ and Alina C.
Iacob$^{3}$
\\{\footnotesize $^{1}$ Department of Mathematics, University of
Kentucky,
Lexington, KY 40506-0027, USA:} 
\\{\footnotesize $^{2}$Departamento de Matem\'atica Aplicada, Universidad de
Murcia, Murcia 30100, Spain:}
   \\{\footnotesize $^{3}$Department
of Mathematics, Georgia Southern University, Statesboro, GA  30460, USA:}
\\{\footnotesize (email:
  enochs@ms.uky.edu; sestrada@um.es;
aiacob@georgiasouthern.edu)}}

\date{}
\begin{document}
\baselineskip=18pt \maketitle
\begin{abstract} 
 Given a double complex $X$  there are spectral sequences with the
$E_2$ terms being either H$_I$ (H$_{II}(X))$ or H$_{II }($H$_I (X))$. But
if $H_I(X)=H_{II}(X)=0$ both spectral sequences have all their terms 0.
This can happen even though there is nonzero (co)homology of interest
associated with $X$. This is frequently the case when dealing with
Tate (co)homology. So in this situation the spectral sequences may not
give any information about the (co)homology of interest. In this article we give a  different way  of constructing homology
groups of $X$  when H$_I(X)=$H$_{II}(X)=0$. With
this result we give a new and elementary proof of balance of Tate
homology and cohomology.

\end{abstract}

\vspace{0.2cm}

{\bf 1.Introduction}\\
 We will mainly be concerned with left $R$-modules over some ring $R$. So
unless otherwise specified, the term module will mean a left $R$-module.
By a complex $(C,d)$ of left $R$-modules we mean a graded module 
$C=(C_n)_{n\in {\bf Z}} $ along with a morphism $d^C=d:C\rightarrow C$ of
graded modules of degree $-1$ such that $d\circ d=0$. We also use
the notation $C=(C^n)$ but where $d$ is of degree $+1$ and where we let
$C_n=C^{-n}$. \\
Given a complex $C$ we let $Z(C)\subset C$ be $Ker(d)$,  let $B(C)
=Im(d)$ and let $H(C)=Z(C)/B(C)$. If $M$ and $N$ are modules and
$C=(C_i)$ and $D=(D^j)$ are complexes, we form complexes denoted $Hom(M,D)$
and $Hom(C, N)$ where $Hom(M,D)^j=Hom(M,D^j)$ and where $Hom(C, N)^i=
Hom(C_i, N)$.\\
By a double complex of modules $X$ we mean a bigraded module $(X^{(i,j)})_{
(i,j)\in {\bf Z}\times {\bf Z}}$ along with morphisms $d'$ and $d''$ of
bidegrees $(1,0)$ and $(0,1)$ respectively such that $d'\circ d'=0$,
$d''\circ d''=0$ and $d'\circ d'' + d''\circ d' =0$.
In [2], $d'$ and $d''$ are denoted $d_1$ and $d_2$ and the homology groups
of $X$ with respect to $d_1$ ($d_2$)  are denoted H$_I(X)$ (H$_{II} (X)$).\\
In this paper we find it convenient to use the related notion of
what Verdier ([5] , Definition 2.1.2) calls a 2-tuple complex and which we will call a
bicomplex. We get a bicomplex if we take the axioms for a double complex
  and  replace the condition $d'\circ d''+d''\circ d'=0$ with 
 the condition that $d'\circ d''= d'' \circ d'$. 
 We can form 
the additive category of bicomplexes where morphisms $f:X\rightarrow Y$
have bidegree $(0,0)$. This category will be an abelian category which
is not only equivalent to but is isomorphic to the category of double
complexes.\\
Given the bicomplex $X$ we let $Z'(X)\subset X$ be the bicomplex $Ker(d')$.
Then we let $B'(X)=Im(d') $ and let $H'(X)$ be the quotient bicomplex
$Z'(X)/B'(X)$. Note that these bicomplexes have their $d'=0$. In a
similar manner we define $Z''(X)$, $B''(X)$ and $H''(X)$. \\
We will tacitly assume that for a bicomplex $X$ we have $X^{(i,j)} \cap
X^{(k,l)}=\emptyset $ unless $(i,j)=(k,l)$. We then will let $x\in X$ mean
that $x\in X^{(i,j)}$  for some (unique) $(i,j)\in {\bf Z}\times {\bf Z}$.
\vspace{.2cm}

{\bf 2. The Main Result}
\vspace{.2cm}

{\bf Theorem 2.1.} {\it Let $X$ be a bicomplex such that $H'(X)=H''(X)=0$ 
i.e. such that $X$ has exact rows and columns.
Then $H'(Z''(X))=H''(Z'(X))$. In this case, if we let $H(X)$ be
$H'(Z''(X))=H''(Z'(X))$, there is a natural isomorphism of bigraded
modules $H(X)\rightarrow H(X)$ of bidegree $(1,-1)$.}
\vspace{.2cm}

{\bf Proof.} It is easy to see that $H'(Z''(X))= Z'(X)\cap Z''(X)/ d'( Z''(X))$
 and that $H''(Z'(X))= Z'(X)\cap Z''(X)/ d''(Z'(X))$. Hence we only need
prove that $d'(Z''(X))=d''(Z'(X))$. We argue this is so by chasing the
diagram. Let $ d'(x)\in d'(Z''(X))$ where $x\in Z''(X)$. Then $d''(x)=0$.
Since $H''(X)=0$ we have $x=d''(y)$ for some $y\in X$.
 Since $d'(d'(y))=0$ we have  $d'(y)\in Z'(X)$. Also $d''(d'(y))=d'(d''(y))
=d'(x)$. So $d'(x)=d'(d''(y))=d''(d'(y))\in d''(Z'(X))$.
 So we have $d' (Z''(X))\subset d''(Z'(X))$. A similar
argument gives that $d''(Z'(X)))\subset d'(Z''(X))$ and so that
$d'(Z''(X))=d''(Z'(X))$.\\
 So when $H'(X)=H''(X)=0$ we let $H(X)=Z'(X)\cap Z''(X)/ d'(Z''(X))= Z'(X)\cap Z''(X)/ d''(Z''(X))$.
In this case we want to find an isomorphism $H(X)\rightarrow H(X)$ of bigraded modules
of bidegree $(1,-1)$. Let $x+d'(Z''(X))\in H(X)=Z'(X)\cap Z''(X)/ d'( Z''(X))$.
Since $x\in Z'(X)\cap Z''(X)$ we have $d''(x)=0$. So since $H''(X)=0$
we have $x=d''(y)$ for some $y\in X$. We claim $d'(y)\in Z'(X)\cap Z''(X)$.
For $d'(d'(y))=0$ and $d''(d'(y))=d'(d''(y))=d'(x)=0$ (since $x\in Z'(X))$.
So we want to map $x+d'(Z''(X))$ to $d'(y)+d'(Z''(X))$. To see that
this map is well-defined, let $\overline{x}+d'(Z''(X))=x+d'(Z''(X))$
where $\overline{x}\in Z'(X)\cap Z''(X)$. Let $d''(\overline{y})=\overline{x}$.
 Then we have $\overline{x}-x
\in d'(Z''(X))$. So let $\overline{x}-x=d'(z)$ where $z\in Z''(X)$.
Then since $H''(X)=0$ we have $w\in X$ with $d''(w)=z$. Then
$d''(d'(w))=d'(d''(w))=d'(z)=\overline{x}-x$. Since $d''(\overline{y}
-y)=\overline{x}-x$ we have $d''(\overline{y}-y -d'(w))=0$, i.e. that
$\overline{y}-y-d'(w) \in Z''(X)$. But $d'(\overline{y}-y-d'(w))=
d'(\overline{y})- d'(y) -0$. So $d'(\overline{y})-d'(y)\in d'(Z''(X))$.
This gives that $x+d'(Z''(X))\mapsto d'(y)+d'(Z''(X))$ (where $d''(y)=x$)
is well-defined. This map is clearly additive, natural and of bidegree $(1,-1)$.
Reversing the roles of $d'$ and $d''$ we get a homomorphism $H(X)
\rightarrow H(X)$ of bidegree $(-1,1)$. By construction we see that
these two maps are inverses of one another, and so both are isomorphisms.

\hfill{$\square$}
\vspace{.2cm}

{\bf 3.  Construction of Bicomplexes}
\vspace{.2cm}

If $C=(C_i)$ and $D=(D^j)$ are complexes we construct a bicomplex
denoted $Hom(C,D)$. We let $Hom(C,D)^{(i,j)}=Hom(C_i, D^j)$ and let
$d'=Hom(d^C, D)$ and $d''=Hom(C, d^D)$. Letting $X=Hom(C,D)$, the 
condition $H'(X)=0$ just says that for each $j\in {\bf Z}$ the complex
$Hom(C, D^j)$ is exact. Similarly the condition $H''(X)=0$ just says
that $Hom(C_i, D)$ is an exact complex for all $i\in {\bf Z}$.\\
 If $X=Hom(C,D)$ 
then want to describe the bicomplexes $Z'(X)$, $Z''(X)$ and $H(X)$ under
certain condition.\\

{\bf Proposition 3.1.} {\it If $C$ is an exact complex and $D$ is any complex, then $Z'(Hom(C,D))
\cong Hom (Z(C), D)$ where the isomorphism is an isomorphism of
bicomplexes.}
\vspace{.2cm}

{\bf Proof.} We have $Z'(Hom(C,D))^{(i,j)}$ is by definition the kernel
of the map $Hom(C_i, D^j)\rightarrow Hom(C_{i+1}, D^j))$. But this
kernel is $Hom(C_i/B_i(C), D^j)$. Since $C$ is exact $C_i/B_i(C)\cong
Z_i(C)$. So as graded modules we have $Z'(Hom(C,D))\cong Hom(Z(C), D)$.
 Clearly these are isomorphisms of bicomplexes.\hfill{$
\square $}
\vspace{.2cm}

We prove the next result in a similar manner.
\vspace{.2cm}

{\bf Proposition 3.2.} {\it If $D$ is an exact complex then $Z''(Hom(C,D))
\cong Hom (C, Z(D))$ where the isomorphism is an isomorphism of bicomplexes.}
\vspace{.2cm}

{\bf Theorem 3.3.} {\it If $C$ and $D$ are both exact complexes and if $Hom(C,D)$
has exact rows and columns then for each $(i,j)\in {\bf Z}\times {\bf Z}$
we have $H(Hom(C,D))^{(i,j)}\cong H^j (Hom(Z_i (C), D))\cong H^i (Hom(C,
Z^j(D)))$.}
\vspace{.2cm}

{\bf Proof.} To get the two isomorphisms we use the two descriptions
of $H(Hom(C,D))$. We first use that $H(Hom(C,D))=H''(Z'Hom(C,D))$.
Since $Z'(Hom(C,D))\cong Hom(Z(C),D)$ as bicomplexes we have
that $H(Hom(C,D))\cong H''(Hom(Z(C),D)$. Since $H''(Hom(Z(C),D))^{(i,j)}=
H^j Hom(Z_i (C), D))$, we get the first isomorphism. Using the other
description of $H(Hom(C,D))$ we get the second isomorphism.\hfill{$\square$}
\vspace{.2cm}

{\bf Corollary 3.4.} {\it For any $n\in {\bf Z}$ we have $H^n (Hom(Z_0(C), D))
\cong H^n (Hom(C, Z^0(D)))$.}
\vspace{.2cm}

{\bf Proof.} Using the natural isomorphism of Theorem 2.1 we have that\\
 $H(Hom(C,D))^{(n,0)}\cong H(Hom(C,D))^{(0,n)}$. 
 Using the two isomorphisms of Theorem  3.3 we have
$H(Hom(C,D))^{(0,n)}\cong H^n (Hom(Z_0(C), D))$ and $H(Hom(C,D))^{(n,0)}
\cong H^n (Hom(C, Z^0 (D)))$. \hfill{$\square$}\\
\vspace{.2cm}

 If $P$  is an exact sequence of projective modules and $E$ an
exact sequence of injective modules, then $Hom(P,E)$ has exact rows
and columns. So this bicomplex satisfies the hypotheses of Theorem 3.3.
 If $M$ and $N$ are a Gorenstein projective and injective
module respectively and if $P$ and $E$ are a complete projective and injective
resolution of $M$ and of $N$ respectively, then  Corollary 3.4 says 
 $H(Hom(M, E))\cong H(Hom(P, N))$. This is what is meant by Tate balance
of cohomology. These (common) cohomology groups are denoted 
$\widehat{Ext}^n (M, N)$ (see [1], section 4 for definitions and notation).

{\bf Remarks 3.5.} The balance over a Gorenstein ring was first proved by Iacob
 ([4], Example 1, pg.2024)
and then by Asadollahi and Salarian [1]
when the ring is local, Gorenstein and the first module finitely generated.
Christensen and Jorgensen in [3] used the inventive
idea of a pinched complex to give a different proof of  the general balance
result. 
\vspace{.2cm}

If we use the tensor product instead of the homomorphism functor, we
get results analogous to the above. The proof of these results will be
easy and obvious modifications of those for the Hom functor. So
we will just state the results. Note that if $C$ and $D$ 
are complexes of left and right $R$-modules respectively, then we
can form a bicomplex which we will denote $C\otimes D$.
\vspace{.2cm}

{\bf Theorem 3.6.} {\it If $C$ is an exact complex of right $R$-modules and
$D$ is an exact complex of left $R$-modules and if the bicomplex $C\otimes
D$ has exact rows and columns then $H^n(Z_0(C) \otimes D)\cong H^n(C\otimes Z_0(D))$ for any $n\in {\bf Z}$.}
\vspace{.2cm}

To get examples where this result can be applied we only need assume
that $C$ and $D$ are both exact complexes of flat modules. Then 
clearly $C\otimes D$ will have exact rows and columns. With this
result we get balance of Tate homology.


\begin{thebibliography}{101}
\bibitem[1]{A1}  Javad Asadollahi and Shokrollah Salarian, {\it Cohomology
theories based on Gorenstein injective modules
}, Trans. Amer. Math. Soc. {\bf 358} (2006), 2183-2203.
\bibitem[2]{2} Henri Cartan and Samuel Eilenberg, Homological Algebra, Princeton
University Press (1956).

\bibitem[3]{A3} Lars Christensen and David Jorgensen, {\it Tate (co)homology
via pinched complexes}, arXiv:1105.2286{\it vi}.

\bibitem[4]{A4} Alina Iacob, {\it Balance in generalized Tate
cohomology}, Comm. Algebra {\bf 33}
(2005), 2009-2024.

\bibitem[5]{A5} Alina Iacob, {\it Absolute, Gorenstein, and Tate torsion 
modules}, Comm. Algebra {\bf 35} (2007), 1589-1606.

\bibitem[6]{A6} Jean-Louis Verdier, {\it Des cat\'{e}gories d\'{e}riv\'{e}es
des cat\'{e}gories ab\'{e}liennes}, Ast\'{e}riques 239 (1997).
\end{thebibliography}
\end{document}